\documentclass[letterpaper, 10 pt, conference]{ieeeconf}  

\IEEEoverridecommandlockouts                              
\usepackage{graphics} 
\usepackage{epsfig} 
\usepackage{times} 
\usepackage{amsmath} 
\usepackage{amssymb}  
\usepackage{algpseudocode}
\usepackage{xcolor}
\usepackage{booktabs}
\usepackage{caption}

\newcommand{\Z}{{\mathbb Z}}

\newcommand{\al}[1]{\begin{align} #1 \end{align}}
\newcommand{\Rs}{\mathbb{R}}
\newcommand{\nn}{\nonumber}

\title{\LARGE \bf Online semi-parametric learning for inverse dynamics modeling}

\author{Diego Romeres$^\dagger$, Mattia Zorzi$^\dagger$,  Raffaello Camoriano$^{\star \ddagger \diamond}$ and Alessandro Chiuso$^\dagger$
\thanks{This work has been  supported by the FIRB project ``Learning meets time'' (RBFR12M3AC) funded by MIUR.}
\thanks{$^\dagger$ Dept. of Information  Engineering, University of Padova, Via Gradenigo 6/b, 35131, Padova, Italy (e-mail: \{\tt \small romeresd,zorzimat,chiuso\}@dei.unipd.it)}%
\thanks{$^{\star}$ iCub Facility, Istituto Italiano di Tecnologia, Via Morego 30, Genoa 16163, Italy (e-mail: \tt \small raffaello.camoriano@iit.it)}
\thanks{$^{\ddagger}$ LCSL, Istituto Italiano di Tecnologia and Massachusetts Institute of Technology, Cambridge, MA 02139, USA.}
\thanks{$^{\diamond}$ DIBRIS, Universit\`a degli Studi di Genova, Via Dodecaneso, 35, Genoa 16145, Italy.}
}

\begin{document}

\maketitle
\thispagestyle{empty}
\pagestyle{empty}

\begin{abstract}
This paper presents a semi-parametric  algorithm for online learning of a robot inverse dynamics model. It combines the strength of the parametric and non-parametric modeling. The former exploits the rigid body dynamics equation, while the latter exploits a suitable kernel function. We provide an extensive comparison with other methods from the literature using real data from the  iCub humanoid robot. In doing so we also compare two different techniques, namely cross validation and marginal likelihood optimization,  for estimating the hyperparameters of the kernel function. 
\end{abstract}


\section{Introduction} 
\label{sec: intro}

Inverse dynamics models are very useful in robotics because they can guarantee high accuracy and low gain control. Building an inverse dynamics model from first principles may be very demanding and, in most cases, out of reach and not suitable for online applications. For this reason, it is of interest to build an inverse dynamics model directly from data, possibly online to allow for real time updating of the model, which is required  for adaptation to changing conditions. 

Traditionally, the inverse dynamics is described by a parametric model given by the rigid body dynamics (RBD), \cite{siciliano2010robotics}. Then, inverse dynamics learning can be recasted as a parametric estimation problem, \cite{hollerbach2008model,Zorzi2014,ZORZI_2015}. The main advantage of this approach is that it provides a global relationship between the input (joint angles, velocities and accelerations) and the output (torques). However, the linear model does not  capture  nonlinearities in the data. To overcome this difficulty, it is possible to describe the inverse dynamics using non-parametric models; we do so by casting the estimation problem in  the Gaussian regression framework, \cite{Rasmussen,BSLCDC,BSL_JOURNAL}, or, equivalently, in the regularization framework,  \cite{rifkin2003regularized}. The latter are characterized by a suitable 
kernel function. However, the drawback of this approach is that a large amount of data is required to produce accurate predictions, as well as high computational load to actually compute the estimated model. This approach is not new and several contributions have recently appeared; for instance  in \cite{ICRA2010NguyenTuong_62320,wu2012semi} the inverse dynamics has been modeled combining the strength of the parametric and of the non-parametric approach. In the latter case two alternatives are possible. The first one is to embed the rigid body dynamics as ``mean'' in the non-parametric part.  The second one is to incorporate the rigid body dynamics in the kernel function.

An important aspect in inverse dynamics learning is the variation of the mechanical properties caused by changing of tasks. It is then necessary to update the model online. In this framework, it is important that the online algorithm is able to take advantage of  the knowledge  already acquired from previously available data, thus speeding up the learning process. This concept is often called transfer learning \cite{pan2010survey, bocsi2013alignment}. Several online learning algorithms have been proposed in the literature. We mention the non-parametric algorithm selecting a sparse subset of training data points (i.e. dictionary), \cite{nguyen2011incremental}, and the semi-parametric algorithms based on the 
locally weighted projection approach, \cite{de2012line}, and on the local Gaussian process regression approach, \cite{nguyen2009model}.
In \cite{gijsberts2011incremental}  a non-parametric online algorithm has been proposed in which the complexity  is kept constant approximating the kernel function using so called ``random features'', \cite{rahimi2007random,quinonero2005unifying}. Finally, in \cite{SEMIPARAMTERIC_2016} a semi-parametric online algorithm, exploiting the above approximation, has been proposed. Here, the rigid body dynamics, preliminarly estimated via least squares, has been embedded as mean in the non-parametric part. 

Another important aspect is the estimation of the hyperparameters of the kernel function. The latter can be estimated according to the maximum likelihood approach, \cite{Rasmussen}, or according to the validation set approach, \cite{james2013introduction}. \\
The first contribution of the paper is to frame various semi-parametric learning  techniques proposed in the literature \cite{ICRA2010NguyenTuong_62320,wu2012semi,SEMIPARAMTERIC_2016} under the same general model, and to provide an online algorithm for this model, exploiting the random features approximation. 
%
%

The second contribution of this paper is to compare these online algorithms 
for estimating the inverse dynamics of right arm of the iCub humanoid robot, \cite{metta2010icub}, \cite{traversaro2013inertial}. In doing that, we also compare the two different approaches for estimating the hyperparameters.

The paper is outlined as follows. In Section \ref{sec: problem_formulation}
we introduce parametric, non-parametric and semi-parametric models.  
In Section \ref{sec: online} the online algorithm to update the model. 
Section \ref{sec:hyperparameter_estimation} deals with the hyperparameters estimation. 
In Section \ref{sec: Simulations} we test the different online algorithms for estimating the inverse dynamics of the right arm of the iCub. 
Finally, in Section \ref{sec:Conclusions} we draw the conclusions.

\section{Inverse Dynamics Learning}
\label{sec: problem_formulation}

Starting from the laws of physics it would in principle be possible to write a (direct) dynamical model which, having as inputs the torques acting on the robot's joints,  outputs the (sampled) trajectory of the free coordinates (joint angles) $q_s$, $s\in \Z$. This is the so called ``direct dynamics''. 

However, for the purpose of control design, it is of interest to know which are the torques  that should be applied in order to obtain a certain trajectory $q_s$. This is the purpose of inverse dynamics modeling: \emph{finding a model which, having joint trajectories as inputs,  outputs the applied torques.}

In order to simplify the modeling exercise, we shall assume, as customary, that not only joint angles $q_s$ can be measured,  but also velocities $\dot{q}_s$ and accelerations $\ddot{q}_s$. Of course this is a (crude) approximation, but we leave possible alternatives to future work. This assumption simplifies considerably the modeling exercise because, given $q_s, \dot{q}_s, \ddot{q}_s$, the inverse dynamics model is, in principle, linear (see \eqref{par_regr_model}).

From now on we shall denote with $x_s=[q_s^\top \; \dot{q}_s^\top \; \ddot{q}_s^\top]^\top\in\Rs^{m}$, $m=3n$,  the vector ``input locations'' 
obtained by stacking positions, velocities and accelerations of all the $n$  joints of the robot. Similarly, $y_s\in\Rs^n$ are the torques applied  to the $n$ joints of the robot at time $s$. 
The inverse dynamics models we consider in this paper will be of the form 
\al{\label{inv:dyn}
 \mathcal{M} :\; \; y_s = h(x_s) +e_s \quad s\in \Z
}
where $h$ is a, possibly non-linear, function and $e_s$ is a zero mean white Gaussian noise with unknown variance $\sigma^2 I_n$. 

\begin{figure}[hbtp]
\centering
\includegraphics[width=\columnwidth]{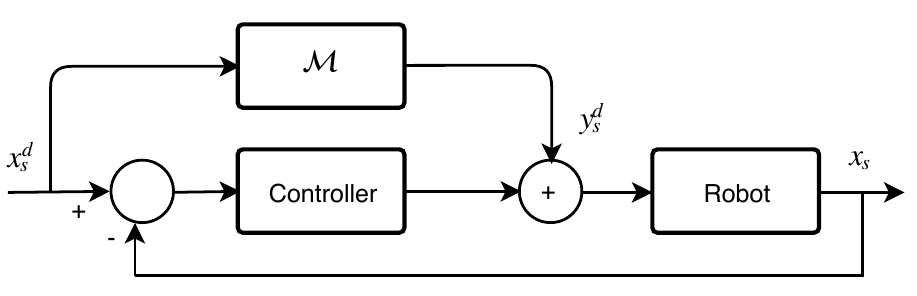}
\caption{Schematic for robot motion control.}\label{fig_aplication}
\end{figure}

The problem of learning the inverse dynamics is that of \emph{estimating} the model $\mathcal{M}$ (i.e. the function $h$) starting from a finite set of measured data samples $\{ y_s,x_s \}_{s=1}^N $.
This model can then be used for robot motion control, see Figure \ref{fig_aplication}.

More precisely, it is exploited to determine the feedforward joint torques $y_s^d$ which should be applied to follow  a desired trajectory $x_s^d$, while employing a feedback controller in order to stabilize the  system. Clearly, the more accurate the inverse dynamics model $\mathcal{M}$ is, the more accurate the motion control is. In this paper we shall consider several approaches, depending upon how the function $h(\cdot)$ in \eqref{inv:dyn} is modelled. 

\subsection{Linear Parametric Model}
The rigid body dynamics (RBD) of a robot is described by the equation
\al{\label{phys_model} y_s= M(q_s) \ddot{q}_s+C(q_s,\dot{q}_s) \ddot{q}_s +G(q_s) }
where $M(q_s)$ is the inertia matrix of the robot, $C(q_s,\dot{q}_s)$ the Coriolis and centripetal forces and $G(q_s)$ the gravity forces, \cite{siciliano2010robotics}. The terms on the right hand side of  (\ref{phys_model}) can be rewritten as $\psi^\top(x_s)\pi$ which is linear in the robot (base) inertial parameters $\pi\in\Rs^{p}$ and where $\psi(x_s)\in\Rs^{p\times n}$ is the known RBD regressor which is a combination of kinematic parameters. In order to make the problem of determining $\pi$ from measured data $y_s$ well posed, we follow a Bayesian approach modeling   $\pi$ as a zero mean Gaussian random vector with covariance matrix $\gamma^2 I_p$. 
Therefore, we consider $h(x) =  \psi(x)^\top \pi$ in \eqref{inv:dyn}, so that 
\al{\label{par_regr_model}y_s= \psi^\top(x_s) \pi+e_s}
where $e_s$ is a zero-mean Gaussian noise with covariance matrix $\sigma^2 I_n$ and it represents nonlinearities of the robot that are not modeled in the rigid body dynamics (e.g. actuator dynamics, friction, etc.). 

\subsection{Nonparametric Model}
The robot inverse dynamics is modeled postulating 
\al{\label{nonpar_model}y_s=g(x_s)+e_s} 
(i.e. $h(x)=g(x)$  in \eqref{inv:dyn}) where 
$g(x_s)$ is a zero mean vector valued (taking values in $\Rs^n$) Gaussian random process  indexed in $\Rs^m$, with covariance function  
\al{\mathbb{E}[g(x_t) g(x_s)^\top]=\rho^2 K(x_t,x_s) I_n.}
 The parameter $\rho^2$ plays the role of  scaling factor and $K$ is a positive definite function, known also as (reproducing) kernel, due to the link between Gaussian process regression and inverse problems in Reproducing Kernel Hilbert Spaces (RKHS) \cite{wahba1990spline}.
In robotics, a typical choice  is the Gaussian kernel, \cite{gijsberts2011incremental,ICRA2010NguyenTuong_62320,wu2012semi},
\al{\label{gaussian_kernel}K(x_t,x_s)=e^{-\frac{\|x_t-x_s\|^2}{2\tau^2}} } 
where $\tau^2$ is the kernel width\footnote{Therefore, to be precise the function $K$ as well as its approximation \eqref{approx_K}, depends on the parameter $\tau$ which will be estimated from data, see Section \ref{sec:hyperparameter_estimation}. For simplicity of exposition this dependence is not made explicit in the notation.}  and it represents the metric to correlate the input locations $x_t$ and $x_s$.
The minimum variance linear estimator of $g$ at time $t$ is given by  the solution of the  regularization problem 
\al{ \label{reg_pb_inf}\hat g_t=\underset{g\in\mathcal{H}}{\mathrm{argmin}} \frac{1}{\sigma^2}\sum_{s=1}^t \|y_s-g(x_s) \|^2+\frac{1}{\rho^2}\| g\|^2_{\mathcal{H}}}
where $\mathcal{H}$ denotes the reproducing kernel Hilbert space (RKHS) 
 of deterministic functions from $\Rs^m$ to $\Rs^n$ associated with $K I_n$ and with norm $\|\cdot \|_{\mathcal{H}}$, \cite{Aronszajn50}.
By the representer Theorem, 
\al{\label{opt_g}\hat g_t(x)=\rho^2\sum_{s=1}^t a_s K(x_s,x) }
where $a _s\in \Rs^n$. Substituting (\ref{opt_g}) in (\ref{reg_pb_inf}) we obtain a Tikhonov regularization problem. 
However, the number of parameters is depending on the number of data $t$, making it hard to obtain on-line (recursive) solutions.  To overcome this limitation, the kernel  $K$ can be approximated, e.g. using the so-called  random features, \cite{rahimi2007random}. This exploits the fact that a positive definite real kernel is the Fourier transform of a non-negative function, which can thus be interpreted as a probability density  \cite{rahimi2007random}. For the Gaussian kernel, that is: 
\al{K(x_t,x_s)= \int_{\Rs^m} p(\omega) e^{i\frac{\omega^\top(x_t-x_s)}{\tau}} \mathrm{d}\omega} 
where \al{\label{pb_density}p(\omega) = \frac{1}{(\sqrt{2\pi} )^m} e^{-\frac{\|\omega\|^2}{2}}.}
 Accordingly, we can approximate $K(x_t,x_s)$ with the sample mean of $e^{i\frac{\omega_k^\top (x_t-x_s)}{\tau}}$, $k=1,..,d$ provided $w_k \sim p(\omega)$, that is: 
 \al{\label{approx_K}K(x_t,x_s)=\frac{1}{d}\sum_{k=1}^d e^{i\frac{\omega_k^\top (x_t-x_s)}{\tau}}=\phi(x_t) ^\top \phi(x_s)} where the basis functions $\phi(x)\in\Rs^{2d}$ are
\al{\phi(x)= & \frac{1}{\sqrt{d}}
\left[\begin{array}{ccc} 
\cos  \left(\frac{\omega_1^\top x}{\tau} \right)  & \ldots &  \cos  \left(\frac{\omega_d^\top x}{\tau} \right)
 \end{array}\right.\nn\\
 &\left.\begin{array}{ccc} 
 \sin \left(\frac{\omega_1^\top x}{\tau} \right)  & \ldots  &  \sin  \left(\frac{\omega_d^\top x}{\tau} \right) 
\end{array}\right]^\top.}
%
%
This is equivalent to model $g(x)$ in the form 
\al{ \label{approx_g} g(x)= (\phi(x)^\top \otimes I_n)  \alpha}
where $\alpha$ is a zero mean Gaussian vector with variance  $\rho^2 I_{2dn}$.
Therefore, the nonparametric model of the robot inverse dynamics \eqref{nonpar_model} can be  approximated by
\al{\label{model_NP} y_s=(\phi(x_s)^\top \otimes I_n ) \alpha+e_s.}
We underline that a peculiarity of model \eqref{model_NP} is that the regressor $\phi(x)$ is depending on the parameter to identify $\tau$.
The advantage of reformulation \eqref{model_NP} is that the dimension of the parameters to estimate $\alpha \in\Rs^{2dn}$ is fixed, which allows a recursive identification formulation, and arbitrary dimensionality: the number of basis functions $d$ can be chosen according to a trade-off between model and computational complexity.  

\subsection{Semi-parametric model with RBD mean} \label{semipar_resid}
This approach combines the parametric and nonparametric models, embedding in the nonparametric model a  mean term, derived from the 
linear parametric model \eqref{par_regr_model}, of the form 
\begin{equation}\label{mean_term}
m_s:= \psi^\top(x_s) \pi
\end{equation}
where 
 $\pi$ is the vector of   inertial parameters and  $\psi(x_s)$ is the RBD regressor. 
Therefore, the inverse dynamics will be modeled as in  \eqref{nonpar_model}  with   $g(x_s)$  a Gaussian process such that 
\begin{equation} \label{SP-Mean} \begin{array}{rl} \mathbb{E}[g(x_s)] &=m_s  = \psi^\top(x_s)\pi \\
\mathrm{Cov}\left( g(x_t), g(x_s) \right) &=\rho^2 K(x_t,x_s)I_n
\end{array}\end{equation}
 where $K$ is the Gaussian Kernel defined in (\ref{gaussian_kernel}). 
  Approximating, as above, the kernel $K$ in \eqref{SP-Mean} with the random features \eqref{approx_K}, the semi-parametric model of the inverse dynamics takes the form
  \al{\label{model_SP} y_s= \psi^\top(x_s) \pi  +(\phi^\top(x_s) \otimes I_n)\alpha+e_s}
 where $\alpha$ is a random vector with zero mean and covariance matrix $\rho^2 I_{2dn}$. 
 As before, $e_s$ is white noise with covariance matrix $\sigma^2 I_n$. 
 
 At this point two  alternatives are possible. The first and most principled one is to treat    $\pi$ as an unknown parameter,  which is to be estimated along with $\rho$, $\sigma$ and $\tau$ using e.g. the marginal likelihood as described in Section \ref{sec:hyperparameter_estimation}. A suboptimal alternative is to assume $\pi$ to be known, possibly estimated using some preliminary experiment as in  \cite{ICRA2010NguyenTuong_62320}.
 In this latter case it will be denoted by $\hat\pi$, and therefore we are only left with modeling the residual vector
  \al{\label{model_SPLS}
 \tilde y_s:=y_s-\psi^\top(x_s) \hat \pi=  (\phi^\top(x_s) \otimes I_n)\alpha+e_s.} 
This latter strategy is followed, for instance, in \cite{SEMIPARAMTERIC_2016}, where the vector $\hat\pi$ is obtained solving in the least squares sense the regression model \eqref{par_regr_model}.

\subsection{Semi-parametric model with RBD kernel} \label{semipar}
An alternative possibility for combining the parametric and nonparametric models in model (\ref{nonpar_model}), is to incorporate the RBD structure in the kernel, \cite{ICRA2010NguyenTuong_62320}. Therefore, $g(x_s)$ is a random process with zero mean and covariance function 
\al{\mathbb{E}[g(x_t)g(x_s)^\top]=\gamma^2 \psi(x_t)^\top \psi(x_s) +\rho^2 K(x_t,x_s) I_n } 
where  the first term $\psi(x_s)$ is the RBD regressor and the second term, $K$, is the Gaussian Kernel defined in (\ref{gaussian_kernel}). As before, $e_s$ is white noise with covariance matrix $\sigma^2 I_n$. Using the kernel approximation (\ref{approx_K}), we have  
\al{\label{cov_SP} \mathbb{E}[g(x_t)g(x_s)^\top]=\gamma^2 \psi(x_t)^\top \psi(x_s) + \rho^2 \phi(x_t)^\top \phi(x_s) I_{n}.}
Accordingly, the approximated semi-parametric model of the inverse dynamics with RBD kernel is:
\al{\label{model_SP2}y_s=\left[\begin{array}{cc} \psi^\top(x_s) & \phi^\top(x_s)\otimes I_n \end{array}\right]\theta+e_s
} where $\theta=[\,\pi^\top \; \alpha^\top\,]^\top \in \Rs^{p+2dn}$ is a zero mean Gaussian random vector  with covariance matrix $\mathrm{blkdiag}(\gamma^2 I_p, \rho^2 I_{2dn})$.  

The semiparametric model with RBD kernel described in this Section is connected, under the Bayesian framework, with the RBD mean model in Section \eqref{semipar_resid}. In fact, model \eqref{model_SP} is equivalent to model \eqref{model_SP2} when $\gamma \rightarrow \infty$ and the parameters $(\rho,\sigma,\tau)$ are fixed; for reasons of space this result will not be explained here. 

\section{Online Learning}
\label{sec: online}

It is apparent that, using the random features approximation of the Gaussian kernel \eqref{approx_K},  all model classes described in the previous Section, see equations \eqref{par_regr_model}, \eqref{model_NP}, \eqref{model_SP} and \eqref{model_SP2},  can ultimately be written in the form  :
\al{ \label{model} \mathcal{M} :\; \; y_s=\varphi(x_s)^\top  \theta+ e_s,\:\: s=1\ldots t}
for a suitable choice of the regressor vector $\varphi(x_s)^\top\in\Rs^{n \times p}$ and $\theta\in\Rs^p$ is modeled as a zero mean random vector with a suitable covariance matrix $\Sigma_0$. $e_s$ is white noise with covariance matrix $\sigma^2 I_n$.  
In this Section we shall assume that $\Sigma_0, \, \tau$ and $\sigma^2$ are known, how to estimate them is a crucial point and will be explained in Section \ref{sec:hyperparameter_estimation}. Thus, the vector $\theta$ completely specifies  the inverse dynamics model and, as such, our learning problem has been reduced to estimating the vector $\theta$ in \eqref{model}.
  At time $t$, the minimum variance linear estimator (i.e. Bayes estimator) of $\theta$ is given by the solution of  the  Tikhonov regularization problem:
\al{\label{ReLS} \hat \theta_t=\underset{\theta\in\Rs^p}{\mathrm{argmin}} \frac{1}{\sigma^2}\sum_{s=1}^t \|y_s-\varphi(x_s)^\top\theta\|^2+\|\theta\|_{\Sigma_0^{-1}}^2.}
This coincides with the so called Regularized Least Squares problem and its optimal solution can be computed recursively through the well known Recursive Least Squares algorithm, see e.g   \cite[Chapter 11]{Ljung}. In practice, the implementation of this algorithm uses Cholesky-based updates \cite{bjorck96},  which have robust numerical properties.

%

\section{Hyperparameters estimation}
\label{sec:hyperparameter_estimation}
All the models presented in Section \ref{sec: problem_formulation} depend on one or more parameter, called hyperparameters, which describe the prior model. For instance, the hyperparameters
in model \eqref{model_SP}, used in semi-parametric learning with RBD mean, are $\xi:=(\pi,\rho^2,\tau^2,\sigma^2)$
while those in model (\ref{model_SP2}), used in semi-parametric learning with RBD kernel, are $\xi:=(\gamma^2,\rho^2,\tau^2,\sigma^2)$.
These hyperparameters are not known  and need to be estimated from the data. In what follows, we consider two different approaches to address this problem. 

\subsection{Validation set approach}
The batch of data used for the identification is split in two data sets: the training set and the validation set. We define a set of candidate hyperparameters and we denote it as $\Xi$. For each $\xi\in\Xi$ we estimate the inverse dynamics model $\mathcal{M}_\xi$ using the training set. Then, for any $\mathcal{M}_\xi$  the mean squared error $\mathrm{MSE}(\xi)$ is computed using the validation set. Hence, the latter provides an estimate of the error rate. 
According to the validation set approach, \cite[Chapter 6]{james2013introduction}, the optimal hyperparameters are given by solving 
\al{\label{CVcost}\hat \xi=\underset{\xi\in\Xi}{\mathrm{argmin}}\; \mathrm{MSE}(\xi).} 

In practice this approach is limited to estimation of a small number of hyperparameters since minimization \eqref{CVcost} is typically performed by gridding the search space $\Xi$.
 
\subsection{Maximum likelihood approach}
Consider, without loss of generality, model (\ref{model}) where both $\theta$ and $e_s$ are assumed to be Gaussian and uncorrelated. Accordingly, the negative marginal loglikelihood (or evidence) of $\mathbf{y}= \left[\begin{array}{ccc} y_1^\top  &\ldots  & y_t^\top \end{array}\right]^\top$ given $\xi\in \Xi$ 
takes the form
\al{L _\xi(\mathbf{y})=\frac{1}{2}\log \det V+ \frac{1}{2} \mathbf{y}^\top V^{-1} \mathbf{y}+c}
where \al{ V=\Phi \Sigma_0\Phi^T+\sigma^2 I_{tn},\;\;\Phi^\top = \left[\begin{array}{ccc} \varphi(x_1) & \ldots &\varphi(x_t)\\  \end{array}\right]\nn}
and $c$ is a term not depending on $\xi$. 
According to the maximum likelihood approach, \cite[Chapter 5]{Rasmussen},
 the optimal hyperparameters are given by solving 
\al{\hat \xi=\underset{\xi\in\Xi}{\mathrm{argmin}}\;L_\xi (\mathbf{y}).} 

\begin{figure}[hbtp]
\centering
\includegraphics[width=\columnwidth]{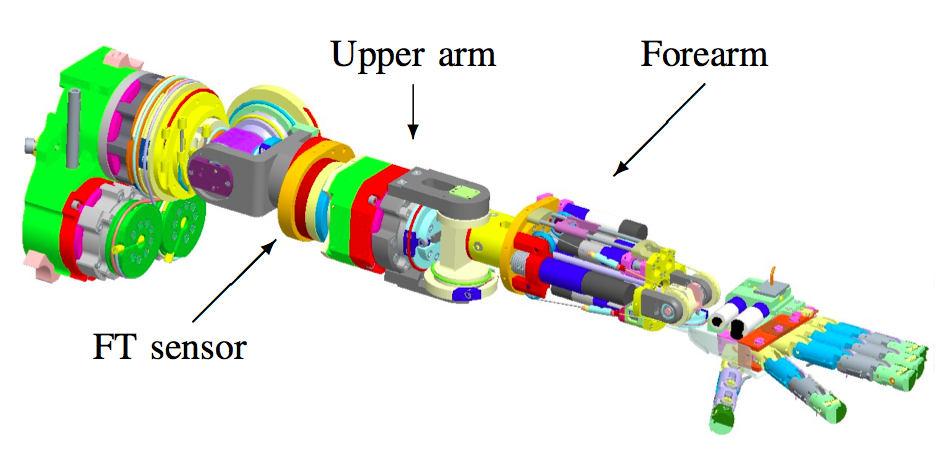}
\caption{iCub's right arm.}\label{fig_right_arm}
\end{figure}

\section{Inverse Dynamics Learning on iCub} \label{sec: Simulations}

iCub is a full-body humanoid robot with 53 degrees of freedom,  \cite{metta2010icub}. We aim to test the models of Section \ref{sec: problem_formulation}  for learning online the inverse dynamics of its right arm.
We consider as inputs $x_s$ the angular positions, velocities and accelerations of the 3 degrees of freedom (dof) shoulder joints and of the 1-dof elbow joint.  
The outputs $y_s$ are the 3 force and 3 torque components measured by the six-axes force/torque (F/T) sensor embedded in the shoulder of the iCub arm, see Figure \ref{fig_right_arm}.

Notice that the measured forces/torques are not the applied joint forces and torques and, as such, the model we learn is not, strictly speaking, the inverse dynamics model. Yet, as explained in \cite{ivaldi2011computing}, the feedforward joint torques can be determined from    components  (forces and  torques) of $y_s$. Indeed, such model has been used in the literature as a benchmark for the inverse dynamics learning, \cite{gijsberts2011incremental}, \cite{SEMIPARAMTERIC_2016}~.

We consider the 2 datasets used in \cite{SEMIPARAMTERIC_2016}, corresponding to different trajectories of the end-effector. In the first one the end-effector tracks circles in the XY  plane  of radius $10 cm$  at an approximative speed of $6m/s$; in the second one, the end-effector tracks similar circles but in the XZ plane (the Z axis corresponds to the vertical direction, parallel to the gravity force). The two circles are tracked using the Cartesian controller proposed in \cite{pattacini2010experimental}. Each dataset contains approximately 8 minutes of data collected  at a sampling rate of $20Hz$, for a total of 10000 points per dataset. One single circle is completed by the robot in about $1.25$ seconds which corresponds to 25 points.

We shall consider  the models described in Section \ref{sec: problem_formulation}, endowed with  the marginal likelihood approach (ML) for the estimation of the hyperparameters, as well as the validation based methods\footnote{As discussed in Section \ref{sec:hyperparameter_estimation}, using validation based methods is unfeasible when the number of hyperparameters is large; therefore we have not applied validation to the semi-parametric model with RBD mean when the mean is to be considered as an hyperparameter nor to the semi-parametric model with RBD kernel which has the extra parameter $\gamma$.}  discussed in \cite{SEMIPARAMTERIC_2016}. For ease of exposition we will use the following shorthands: 
\begin{itemize}
\item
\emph{P}:  the parametric model.
\item \emph{NP-ML}:  the nonparametric model; hyperparameters estimated  with ML.
\item \emph{SP-ML}: the semi-parametric model with RBD mean;  hyperparameters estimated  with ML.
\item \emph{SP2-ML}:  the semi-parametric model with RBD mean, in which the mean  is computed via least squares as in \cite{SEMIPARAMTERIC_2016} and then the nonparametric model is applied to the residuals (see Section \eqref{semipar_resid}); hyperparameters estimated  with ML.
\item \emph{SPK-ML}: the semi-parametric model with RBD kernel; hyperparameters estimated  with ML.
\item \emph{NP-VS}: the nonparametric model with hyperparameters estimated with VS.
\item \emph{SP2-VS}: the semi-parametric model with RBD mean, in which the mean  is computed via least squares; hyperparameters estimated with VS.

\end{itemize}

The proposed algorithms have been implemented using Matlab. The RBD regressor $\psi$ for the right arm of iCub has been computed using the  library iDynTree, \cite{nori2015icub}. The Marginal Likelihood has been optimized using the Matlab \verb!fminsearch.m! function.
The recursive least square algorithms have been implemented using GURLS library, \cite{tacchetti2013gurls}. 
The results of all validation based methods are obtained using  code which has been  kindly provided by the authors of \cite{SEMIPARAMTERIC_2016}.

For each algorithm as above,  we consider the following online learning scenario (with reference to the general model structure \eqref{model}):
\begin{itemize}
\item Initialization: The first 1000 points in XY-dataset are used to estimate the hyperparameters, as well as to compute an    
 initial estimate of parameter $\theta$, say $\hat{\theta}_0$.
 \item Training XY:  Use the remaining 9000 points of XY-dataset to update online parameter $\theta$ using the recursive least-squares algorithm, thus obtaining   $\hat \theta_t$, $t=1,\dots,9000$.
\item Training XZ: The  XZ-dataset is split in 5 sequential subsets of 2000 points (approximately $100$ seconds) each. For each subset we update online the parameter $\theta$ independently always initializing the recursions with  $\hat \theta_{9000}$, computed from the  training dataset XY .
\end{itemize}

In the last step of the procedure,  the initial model has been computed from the Training XY dataset, which corresponds to a different motion with respect to   XZ-dataset. Our goal is that the model estimated with the second dataset quickly captures the new information gathered from the XZ-dataset,  adapting to the new task. For instance, in model predictive control the quality of the control 
depends on the prediction capability of the model over a prescribed horizon, \cite{maciejowski2002predictive}. In order to measure this ability we consider the following index:
\begin{align}
\label{eq:pred_error}
\varepsilon^{(k)}_t &= {\frac{\sum_{s=1}^{T}(y^{(k)}_{t+s} -\hat{y}^{(k)}_{t+s|t} )^2}{\sum_{s=1}^{T} (y^{(k)}_{t+s})^2}}
\end{align}
where $\hat{y}_{t+s|t}^{(k)}$ is the estimate of  the output  $y^{(k)}_{t+s}$ at time $t+s$ using the \emph{model} estimated with data up to time $t$.   
Therefore, $\varepsilon^{(k)}_t$ represents the relative squared prediction error over the horizon $[t+1,\ldots,t+T]$ using model $\mathcal{M}$
computed at time $t$. Let $\varepsilon_t^F$ and $\varepsilon_t^T$ be the average value of $\varepsilon_t^{(k)}$ for the 3 forces and the 3 torques, respectively. \begin{figure}[hbtp]
\centering
\includegraphics[scale=0.5]{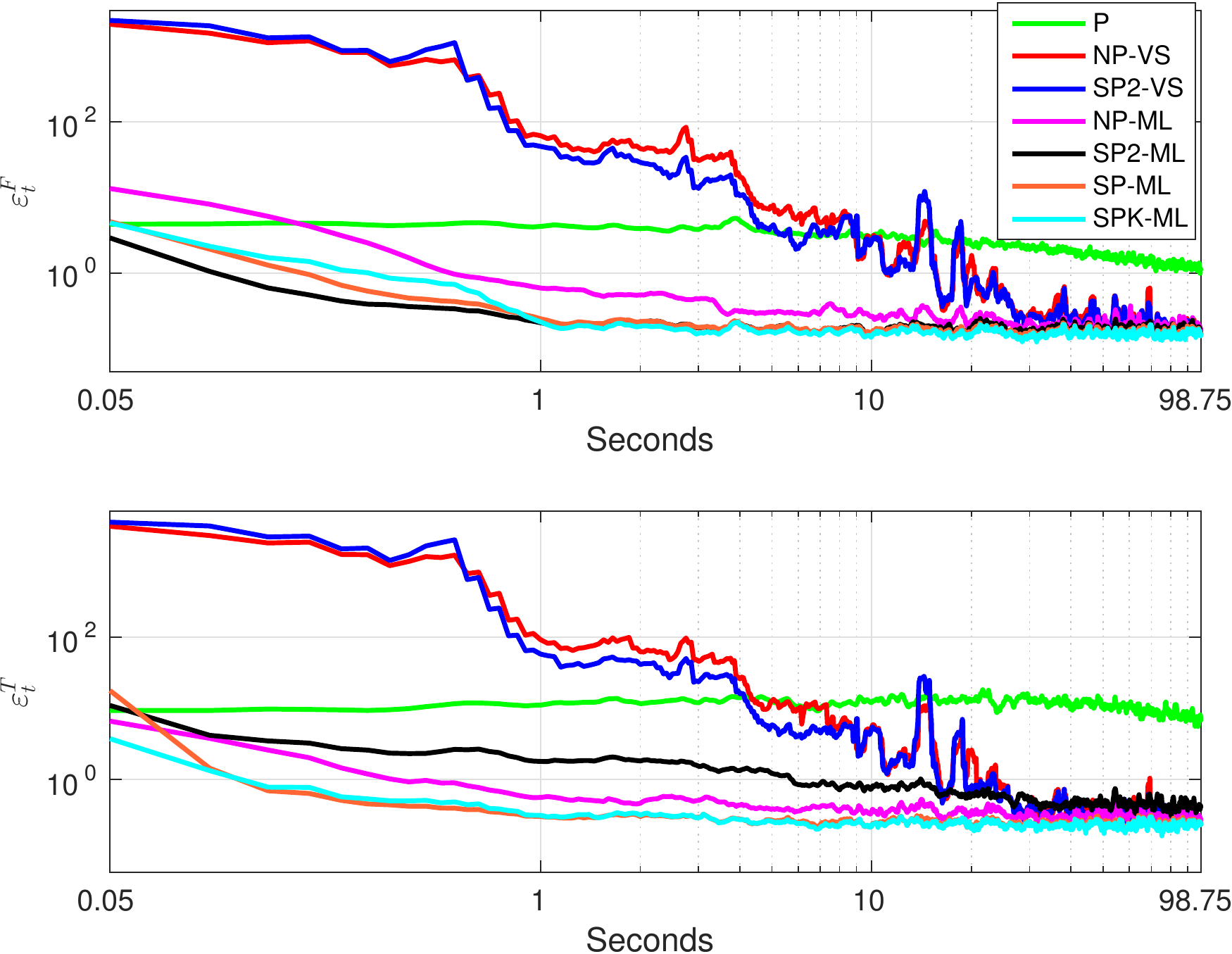}
\caption{Average (over the 5 subsets of 100 seconds each) of the relative squared prediction errors $\varepsilon_t^F$ and $\varepsilon_t^T$, computed with $T=25$ corresponding to an  horizon of 1.25 seconds.}
\label{fig:prediction_error}
\end{figure}

\begin{figure}[hbtp]
\centering
\includegraphics[scale=0.5]{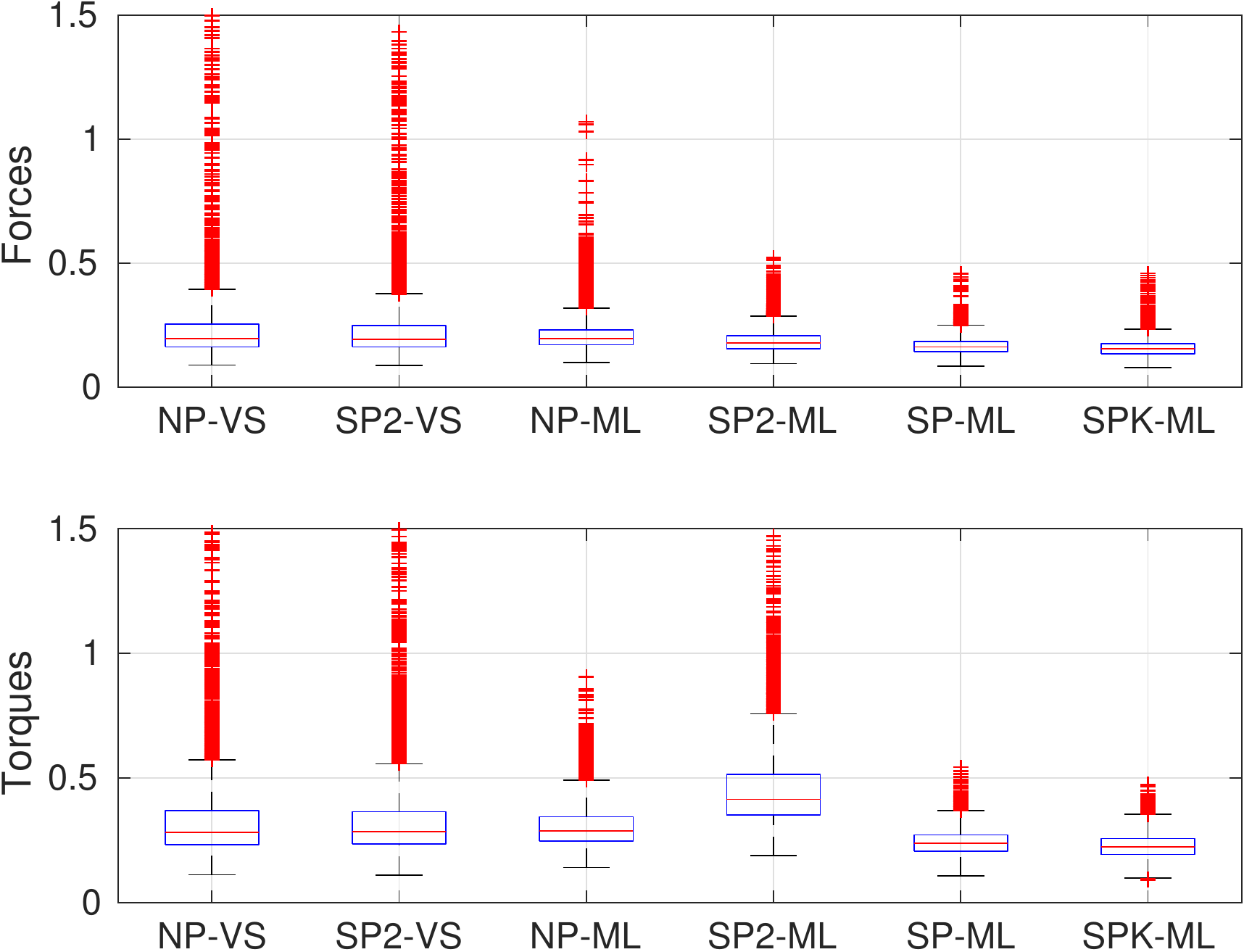}
\caption{Boxplots of the steady state (i.e.  after $30$ seconds, see Figure \ref{fig:prediction_error}) relative squared prediction errors 
$\varepsilon_t^F$ and $\varepsilon_t^T$, computed with $T=25$ corresponding to an  horizon of 1.25 seconds.}
\label{fig:prediction_error_boxplot}
\end{figure}

In Figure \ref{fig:prediction_error} we show $\varepsilon_t^F$ and $\varepsilon_t^T$, averaged over the $5$ subsets, with $T=25$ (1.25 seconds), i.e. with  the end-effector completing one circle during the prediction horizon. Clearly, the parametric algorithm P exhibits a poor performance because it describes only crude idealizations of the actual dynamics. The algorithms based on the VS approach perform significantly worse in the first $60$ seconds than those based on the ML approach. This result is not unexpected because the ML approach represents a robust way to estimate hyperparamters, \cite{Pillonetto2015106}. The models with the best performance are SP-ML and SPK-ML because they combine the benefit of the parametric and the non-parametric approach. Although also SP2-ML exploits this benefit, it provides a slightly worse performance. This is probably due by the fact that  the first (least squares) step, i.e. estimation of the linear model, is subject to a strong bias deriving from the unmodeled dynamics. Instead, a sound approach is followed by SP-ML and SPK-ML  in which the estimation of the hyperparameters is performed  jointly, avoiding such bias. In the steady state all these methods, with the exception of P, provide similar performance; yet the two semi-parametric models (SP-ML and SPK-ML) perform better both in terms of average as well as standard deviation, as clearly shown in  Figure \ref{fig:prediction_error_boxplot}
which reports the boxplots of $\varepsilon_t^F$ and $\varepsilon_t^T$ in ``steady state'', i.e. after the first $30$ seconds which is considered to be transient (see Figure  \ref{fig:prediction_error}).


\section{Conclusions} \label{sec:Conclusions}
 In this paper we have placed several algorithms used for online learning of the robot inverse dynamics in a common framework. Such algorithms are classified according to the considered model (parametric, non-parametric, semi-parametric with RBD mean and semi-parametric with RBD kernel) and according to the way the hyperparameters are estimated (VS approach and ML approach).      
 We applied those algorithms for online leaning of the inverse dynamics of the right arm of the iCub. 
 The results showed the superiority of the ML approach to estimate the hyperparameters. Finally, semi-parametric models outperform the others.  
 The latter result confirms the advantage in 
combining parametric and non-parametric approaches together.  
\section{ACKNOWLEDGMENTS}
The authors gratefully acknowledge the iCub Facility and LCSL-IIT@MIT research groups led by Francesco Nori and Lorenzo Rosasco for making their data and code available to us.




\begin{thebibliography}{10}
\providecommand{\url}[1]{#1}
\csname url@samestyle\endcsname
\providecommand{\newblock}{\relax}
\providecommand{\bibinfo}[2]{#2}
\providecommand{\BIBentrySTDinterwordspacing}{\spaceskip=0pt\relax}
\providecommand{\BIBentryALTinterwordstretchfactor}{4}
\providecommand{\BIBentryALTinterwordspacing}{\spaceskip=\fontdimen2\font plus
\BIBentryALTinterwordstretchfactor\fontdimen3\font minus
  \fontdimen4\font\relax}
\providecommand{\BIBforeignlanguage}[2]{{%
\expandafter\ifx\csname l@#1\endcsname\relax
\typeout{** WARNING: IEEEtran.bst: No hyphenation pattern has been}%
\typeout{** loaded for the language `#1'. Using the pattern for}%
\typeout{** the default language instead.}%
\else
\language=\csname l@#1\endcsname
\fi
#2}}
\providecommand{\BIBdecl}{\relax}
\BIBdecl

\bibitem{siciliano2010robotics}
B.~Siciliano, L.~Sciavicco, L.~Villani, and G.~Oriolo, \emph{Robotics:
  modelling, planning and control}.\hskip 1em plus 0.5em minus 0.4em\relax
  Springer Science \& Business Media, 2010.

\bibitem{hollerbach2008model}
J.~Hollerbach, W.~Khalil, and M.~Gautier, ``Model identification,'' in
  \emph{Springer Handbook of Robotics}.\hskip 1em plus 0.5em minus 0.4em\relax
  Springer, 2008, pp. 321--344.

\bibitem{Zorzi2014}
M.~Zorzi, ``Rational approximations of spectral densities based on the alpha
  divergence,'' \emph{Mathematics of Control, Signals, and Systems}, vol.~26,
  no.~2, pp. 259--278, 2014.

\bibitem{ZORZI_2015}
------, ``{M}ultivariate {S}pectral {E}stimation based on the concept of
  {O}ptimal {P}rediction,'' \emph{IEEE Trans. Autom. Control}, vol.~60, pp.
  1647--1652, Jun. 2015.

\bibitem{Rasmussen}
C.~Rasmussen and C.~Williams, \emph{{G}aussian Processes for Machine
  Learning}.\hskip 1em plus 0.5em minus 0.4em\relax The MIT Press, 2006.

\bibitem{BSLCDC}
M.~Zorzi and A.~Chiuso, ``A {B}ayesian approach to sparse plus low rank network
  identification,'' in \emph{54th IEEE Conference on Decision and Control}, Dec
  2015, pp. 7386--7391.

\bibitem{BSL_JOURNAL}
------, ``{S}parse plus {L}ow rank {N}etwork {I}dentification: A
  {N}onparametric {A}pproach,'' \emph{Automatica}, vol.~53, 2017.

\bibitem{rifkin2003regularized}
R.~Rifkin, G.~Yeo, and T.~Poggio, ``Regularized least-squares classification,''
  \emph{Nato Science Series Sub Series III Computer and Systems Sciences}, vol.
  190, pp. 131--154, 2003.

\bibitem{ICRA2010NguyenTuong_62320}
D.~Nguyen-Tuong and J.~Peters, ``Using model knowledge for learning inverse
  dynamics,'' in \emph{IEEE International Conference on Robotics and
  Automation}, 2010.

\bibitem{wu2012semi}
T.~Wu and J.~Movellan, ``Semi-parametric gaussian process for robot system
  identification,'' in \emph{IEEE/RSJ International Conference on Intelligent
  Robots and Systems (IROS)}, 2012, pp. 725--731.

\bibitem{pan2010survey}
S.~J. Pan and Q.~Yang, ``A survey on transfer learning,'' \emph{IEEE
  Transactions on Knowledge and Data Engineering}, vol.~22, no.~10, pp.
  1345--1359, 2010.

\bibitem{bocsi2013alignment}
B.~Bocsi, L.~Csat{\'o}, and J.~Peters, ``Alignment-based transfer learning for
  robot models,'' in \emph{The 2013 International Joint Conference on Neural
  Networks (IJCNN)}, 2013, pp. 1--7.

\bibitem{nguyen2011incremental}
D.~Nguyen-Tuong and J.~Peters, ``Incremental online sparsification for model
  learning in real-time robot control,'' \emph{Neurocomputing}, vol.~74,
  no.~11, pp. 1859--1867, 2011.

\bibitem{de2012line}
J.~S. de~la Cruz, D.~Kulic, W.~S. Owen, E.~Calisgan, and E.~A. Croft, ``On-line
  dynamic model learning for manipulator control.'' in \emph{SyRoCo}, 2012, pp.
  869--874.

\bibitem{nguyen2009model}
D.~Nguyen-Tuong, M.~Seeger, and J.~Peters, ``Model learning with local gaussian
  process regression,'' \emph{Advanced Robotics}, vol.~23, no.~15, pp.
  2015--2034, 2009.

\bibitem{gijsberts2011incremental}
A.~Gijsberts and G.~Metta, ``Incremental learning of robot dynamics using
  random features,'' in \emph{IEEE International Conference on Robotics and
  Automation (ICRA)}, 2011, pp. 951--956.

\bibitem{rahimi2007random}
A.~Rahimi and B.~Recht, ``Random features for large-scale kernel machines,'' in
  \emph{Advances in neural information processing systems}, 2007, pp.
  1177--1184.

\bibitem{quinonero2005unifying}
J.~Quinonero-Candela and C.~E. Rasmussen, ``A unifying view of sparse
  approximate gaussian process regression,'' \emph{The Journal of Machine
  Learning Research}, vol.~6, pp. 1939--1959, 2005.

\bibitem{SEMIPARAMTERIC_2016}
R.~Camoriano, S.~Traversaro, L.~Rosasco, G.~Metta, and F.~Nori, ``Incremental
  semiparametric inverse dynamics learning,'' in \emph{2016 IEEE International
  Conference on Robotics and Automation (ICRA)}, May 2016, pp. 544--550.

\bibitem{james2013introduction}
G.~James, D.~Witten, T.~Hastie, and R.~Tibshirani, \emph{An introduction to
  statistical learning}.\hskip 1em plus 0.5em minus 0.4em\relax Springer, 2013,
  vol. 112.

\bibitem{metta2010icub}
G.~Metta, L.~Natale, F.~Nori, G.~Sandini, D.~Vernon, L.~Fadiga, C.~Von~Hofsten,
  K.~Rosander, M.~Lopes, J.~Santos-Victor \emph{et~al.}, ``The icub humanoid
  robot: An open-systems platform for research in cognitive development,''
  \emph{Neural Networks}, vol.~23, no.~8, pp. 1125--1134, 2010.

\bibitem{traversaro2013inertial}
S.~Traversaro, A.~Del~Prete, R.~Muradore, L.~Natale, and F.~Nori, ``Inertial
  parameter identification including friction and motor dynamics,'' in
  \emph{13th IEEE-RAS International Conference on Humanoid Robots (Humanoids)},
  2013, pp. 68--73.

\bibitem{wahba1990spline}
G.~Wahba, \emph{Spline models for observational data}.\hskip 1em plus 0.5em
  minus 0.4em\relax Siam, 1990, vol.~59.

\bibitem{Aronszajn50}
N.~Aronszajn, ``Theory of reproducing kernels,'' \emph{Trans. of the American
  Mathematical Society}, vol.~68, pp. 337--404, 1950.

\bibitem{Ljung}
L.~Ljung, \emph{System Identification, Theory for the User}.\hskip 1em plus
  0.5em minus 0.4em\relax Prentice Hall, 1997.

\bibitem{bjorck96}
\BIBentryALTinterwordspacing
A.~Bj\"orck, \emph{Numerical Methods for Least Squares Problems}.\hskip 1em
  plus 0.5em minus 0.4em\relax Society for Industrial and Applied Mathematics,
  1996. [Online]. Available:
  \url{http://epubs.siam.org/doi/abs/10.1137/1.9781611971484}
\BIBentrySTDinterwordspacing

\bibitem{ivaldi2011computing}
S.~Ivaldi, M.~Fumagalli, M.~Randazzo, F.~Nori, G.~Metta, and G.~Sandini,
  ``Computing robot internal/external wrenches by means of inertial, tactile
  and f/t sensors: theory and implementation on the icub,'' in \emph{Humanoid
  Robots (Humanoids), 2011 11th IEEE-RAS International Conference on}, 2011,
  pp. 521--528.

\bibitem{pattacini2010experimental}
U.~Pattacini, F.~Nori, L.~Natale, G.~Metta, and G.~Sandini, ``An experimental
  evaluation of a novel minimum-jerk cartesian controller for humanoid
  robots,'' \emph{IEEE/RSJ International Conference on Intelligent Robots and
  Systems}, pp. 1668--1674, 2010.

\bibitem{nori2015icub}
F.~Nori, S.~Traversaro, J.~Eljaik, F.~Romano, A.~Del~Prete, and D.~Pucci,
  ``icub whole-body control through force regulation on rigid non-coplanar
  contacts,'' \emph{Frontiers in Robotics and AI}, p.~18, 2015.

\bibitem{tacchetti2013gurls}
A.~Tacchetti, P.~Mallapragada, M.~Santoro, and R.~Rosasco, ``Gurls: A least
  squares library for supervised learning,'' \emph{Journal of Machine Learning
  Research}, vol.~14, pp. 3201--3205, 2013.

\bibitem{maciejowski2002predictive}
J.~M. Maciejowski, \emph{Predictive control: with constraints}.\hskip 1em plus
  0.5em minus 0.4em\relax Pearson education, 2002.

\bibitem{Pillonetto2015106}
G.~Pillonetto and A.~Chiuso, ``Tuning complexity in regularized kernel-based
  regression and linear system identification: The robustness of the marginal
  likelihood estimator,'' \emph{Automatica}, vol.~58, pp. 106 -- 117, 2015.

\end{thebibliography}
\end{document}